\documentstyle[twoside,amssymb,amsmath,12pt]{article}
\def\C{{\Bbb C}}

\def\R{{\Bbb R}}
\def\P{{\Bbb P}}
\def\Z{{\Bbb Z}}
\def\A{{\Bbb A}}

\newtheorem{prop}{Proposition}[section]
\newtheorem{dfn}[prop]{Definition}
\newtheorem{theo}[prop]{Theorem}

\newtheorem{rem}[prop]{Remark}

\newtheorem{exam}[prop]{Example}
\newtheorem{ques}[prop]{Question}
\newtheorem{prob}[prop]{Problem}

\title{\sc Toric Degenerations of Fano Varieties and Constructing 
Mirror Manifolds} 

\author{{\sc Victor V. Batyrev} \\
\small  {\em Mathematisches Institut, Universit\"at T\"ubingen}   \\
\small  {\em Auf der Morgenstelle 10,  72076  T\"ubingen, Germany}  \\
\small  {\em e-mail: batyrev@bastau.mathematik.uni-tuebingen.de} \\
 }

\begin{document}

\date{}

\maketitle

\begin{abstract}
For an arbitrary  smooth $n$-dimensional Fano variety $X$ 
we introduce  the notion of a {\em small toric degeneration}. 
Using small toric degenerations of Fano $n$-folds $X$,  we propose 
a general method for  constructing  mirrors of  Calabi-Yau complete 
intersections in $X$. Our mirror construction  is based on a generalized  
monomial-divisor mirror correspondence  which can be used  for 
computing Gromov-Witten invariants  of rational curves via specializations
of GKZ-hypergeometric series\footnote{This is an extended version of the 
author's talk given during the  Summer Symposium 
on Algebra at University of Niigata, July 22-25, 1997.}. 
\end{abstract}

\section{Introduction} 

Recent progress in understanding the mirror symmetry phenomenon 
using explicit mirror constructions for Calabi-Yau hypersurfaces 
and complete intersections in toric varieties \cite{BA,BS,BB0,Bo}  
leads to  the following natural question: 
\medskip

{\em Is it possible to extend the mirror constructions  for Calabi-Yau 
complete intersections in toric Fano varieties to 
the case of Calabi-Yau complete intersections in nontoric Fano varieties? }
\medskip 

The first progress  in this direction 
has been obtained for Grassmannians \cite{BCKS1} and,  more generally, 
for partial flag manifolds \cite{BCKS2}. The key idea in both 
examples is based on a degeneration of Grassmannians (resp.  
partial flag manifolds) to some singular Gorenstein 
toric Fano varieties. These
degenerations have been introduced and investigated by Sturmfels, 
Gonciulea and Lakshmibai  in \cite{GL1,GL2,L,S1,S2}. 

The present paper is aimed to give a short systematic overview of our 
method for   constructing mirror manifolds and to formulate 
some naturally arising  questions and open problems.

In Section 2 we start with a review of a method for 
constructing degenerations of unirational varieties $X$ to toric 
varieties $Y$ using canonical subalgebra bases. This method has been 
discovered  by Kapur \& Madlener \cite{KM} and independently 
by Robbiano \& Sweedler \cite{RS}. Further results on this topic have 
been obtained  in \cite{O,M,S1} (see also \cite{S2} for more details).  

In Section 3 we introduce the notion of a 
{\em small toric degeneration} of a Fano manifold and 
discuss some examples. Finally, in Section 4 we explain our generalized 
mirror construction which uses small toric degenerations.

\section{Canonical subalgebra bases}

Let $A$ be a finitely generated subalgebra of the polynomial 
ring 
$$K[{\bf u}]:=K[u_1, \ldots, u_n],$$
i.e., $X= Spec\, A$ is an  unirational  affine algebraic  variety 
together with  a dominant morhism  $\A^n \to X$.   
We choose  a weight vector 
$\omega = (\omega_1, \ldots, \omega_d) \in \R^n$ and  
set $$wt({\bf u}^{\bf a}) = wt(u_1^{a_1} \cdots u_n^{a_n}) := 
 \sum_{i=1}^n a_i \omega_i.$$
The number $wt({\bf u}^{\bf a})$ will be called the {\bf weight} of the
monomial  ${\bf u}^{\bf a}$.
We define  
a partial order on the set of all monomials in $K[{\bf u}]$ as follows: 
\[ {\bf u}^{\bf a} \prec  {\bf u}^{\bf a'} \Leftrightarrow 
wt({\bf u}^{\bf a}) \leq   wt({\bf u}^{\bf a'}). \]
If $f \in K[{\bf u}]$ is a polynomial, then $in_{\prec}(f)$ 
denotes the {\bf initial part of $f$}, i.e., the sum of 
those monomials in $f$ whose weight is maximal. 
By definition, one has $in_{\prec}(fg) = in_{\prec}(f) in_{\prec}(g)$. 
For suficiently general choice of the weight vector $\omega \in \R^n$ 
the initial part of a polynomial  $f \in K[{\bf u}]$ is a single 
monomial.

\begin{dfn} 
{\rm The $K$-vector space spanned by initial terms of elements $f \in A$ 
is called the {\bf  initial algebra} and is denoted by 
\[ in_{\prec}(A) : = \{ in_{\prec}(f)\; : \; f \in A \}. \]}
\end{dfn}

\begin{dfn} 
{\rm A subset ${\cal F} \subset A $ is called a {\bf canonical basis of 
the subalgebra} $A \subset K[{\bf u}]$, if the initial 
subalgebra $in_{\prec}(A)$ is generated 
by the elements  
$$\{ in_{\prec}(f)\; : \; f \in {\cal F} \}.$$ } 
\end{dfn}

Fix a set of 
polynomials ${\cal F} = \{f_1, \ldots, f_m\} \subset A$.  
We set $K[{\bf v}]:=K[v_1, \ldots, v_m]$. 
Let $I$ be the kernel of the 
canonical epimorphism 
$$\varphi \; : \; K[{\bf v}] \to A$$ 
$$v_i \mapsto  f_i$$ and $I_{\prec}$ 
the kernel of the 
canonical epimorphism 
$$\varphi_0 \; : \; K[{\bf v}]  \to in_{\prec}(A)$$ 
$$v_i \mapsto in_{\prec}(f_i)$$

\begin{rem} 
{\rm It is easy to show that the ideal $I_{\prec}$ is generated by 
binomials (see \cite{ES} for 
general theory of binomial ideals). Hence, the spectrum 
of $in_{\prec}(A)$ is 
an affine toric variety (possibly not normal).} 
\end{rem}  
 
Now we assume that  $\omega = (\omega_1, \ldots, \omega_d) \in \Z^n$ an integral  
weight vector. If the set of polynomials 
${\cal F} = \{f_1, \ldots, f_m\} \subset A$ form a canonical 
basis of the subalgebra $A \subset K[{\bf u}]$ with respect to the 
partial order defined by  $\omega$, then we can define a $1$-parameter 
family of subalgebras 
\[ A_t := \{ f(t^{-\omega_1}u_1, \ldots, t^{-\omega_n} u_n) \;\; |\; \; 
f(u_1, \ldots, u_n) \in A \},\; \;\; t \in K \setminus \{ 0 \} \}. \] 
Setting $A_0: = in_{\prec}(A)$, we   obtain a flat 
family of subalgebras $A_t \subset  K[{\bf u}]$ such that 
$A_t \cong A$ for $ t\neq 0$ and $A_0 \cong  K[{\bf v}]/I_{\prec}$. 
This allows us to consider the affine   toric 
variety  $Spec\, A_0$ as a flat degeneration of 
$Spec\, A$.

\begin{rem} 
{\rm It is important to remark that the above method for constructing 
toric degenerations strongly depends on the choice of the coordinates 
$u_1, \ldots, u_n$ on $\A^n$ and on the choice of a weight vector $\omega$.} 
\end{rem}

\begin{exam} 
{\rm 
Let $A(r,s) \subset K[{\bf X}]: = K[X_{ij}]$ $(1 \leq i \leq r, \; 
1 \leq j \leq s)$ be the subalgebra of the polynomial algebra 
$K[{\bf X}]$ generated by all $r \times r$ minors of a generic 
$r \times s$ matrix $( r \leq  s)$, i.e.,  $A(r,s)$ is the homogeneous 
cooordinate ring of the Pl\"ucker embedded Grassmannian 
$G(r,s) \subset \P^{ { s \choose r } -1}$. Define the weights of monomials 
as follows
\[ wt(X_{ij}) := (j-1) s^{i-1}, \;\; i,j \geq 1.\]
In particular, one has 
\[ wt(X_{1,i_1} \cdots X_{r,i_r}) = (i_1-1) + (i_2-1)s + \cdots + 
(i_r-1)s^{r-1} \]
and therefore the initial term of each 
 $(i_1, \ldots, i_r)$-minor $(1 \leq i_1 < \cdots < i_r \leq s)$ is 
exactly the product of terms on the main diagonal: 
\[X_{1,i_1} \cdots X_{r,i_r}. \]  

The following result  is due to Sturmfels \cite{S1,S2}:

\begin{theo}
The set of all  $s \times s$-minors form a canonical base of  the 
subalgebra  $A(r,s) \subset K[{\bf X}]$ with respect to the 
partial order defined by the above weight vector. In particuar, one 
obtains a natural toric degeneration of the Grassmanninan $G(r,s)$.  
\end{theo} 

\label{grass}
} 
\end{exam} 

\section{Small toric degenerations of Fano varieties} 

\begin{dfn} 
{\rm Let $X \subset \P^m$ be a smooth Fano variety of dimension $n$. A normal 
Gorenstein 
toric Fano variety $Y \subset \P^m$ is called a 
{\bf small toric degeneration} of $X$, if 
there exists a Zariski open neighbourhood  $U$ of $0  \in  \A^1$  and 
an irreducible subvariety ${\frak X} \subset \P^m \times U$ such that 
the   morphism   $\pi\; : \; {\frak X} \to U$ is flat and 
the following conditions hold: 

{(i)} the fiber $X_t := \pi^{-1}(t) \subset \P^m$ is smooth 
for all $t \in U \setminus \{ 0 \}$; 

{(ii)} the special fiber $X_0 := \pi^{-1}(0) \subset \P^m$ has  
at worst Gorenstein terminal singularities (see \cite{KMM}) 
and $X_0$ is isomorphic to $Y \subset \P^m$; 

{(iii)} the canonical homomorphism 
\[ Pic({\frak X}/U) \to Pic(X_t) \]
is an isomorphism for all $t \in U$. } 
\label{def-small}
\end{dfn} 

\begin{rem} 
{\rm It is weill-known that if $Y$ has at worst terminal singularities, then 
the codimension of the singular locus of $Y$ is at least $3$. 
On the other hand, it is easy to show that the only possible toric 
Gorenstein terminal singularities in dimension $3$ 
are ordinary double points (or nodes): $x_1x_2 - x_3x_4=0$. So, if $Y$ 
is a small toric degeneration of $X$, then the singular locus of $Y$ in 
codimension $3$ must consist of nodes.}
\label{codim3}
\end{rem} 

\begin{exam} 
{\rm Let $Y:= P(r,s) \subset \P^{ { s \choose r } -1}$ be 
the toric degeneration of the Grassmannian $X:= Gr(r,s) \subset \P^{ { s \choose r } -1}$ (see Example \ref{grass}). Then  
$Y$ is a small toric degeneration of $X$ \cite{BCKS1}. }
\end{exam}     

\begin{exam} 
{\rm Let $X:= F(n_1, \ldots,n_k ,n) \subset \P^{m}$ be the partial flag
manifold it is Pl\"ucker embedding. It is proved in 
\cite{BCKS2} that the toric degenerations 
introduced and investigated by 
Gonciulea and Lakshmibai  in \cite{GL1,GL2,L} are small toric degenerations 
of $X$.
} 
\end{exam}

\begin{exam} 
{\rm Let $V_{d,n} \subset \P^{n+1}$ 
be a Gorenstein toric Fano hypersurface of degree $d$ 
$(d \geq 2)$ 
in projective space of dimension $n \geq 2d -2$  defined by 
the homogeneous equation 
\[ z_1 \cdots z_d = z_{d+1} \cdots z_{2d}. \]
It is easy to check that 
 irreducible components of the singular locus of  $V_{d,n}$
 are  
\[ \frac{d^2(d-1)^2}{4} \]
codimension-3 linear subspaces 
\[ z_i=z_j=z_k =z_l =0, \]
\[ \; \; \{i,j \} \subset \{1, \ldots, d\},\; 
 \{k,l \} \subset \{d+1, \ldots, 2d \}, \;i \neq j, \; k \neq l. \]
consisting of  nodes.} 
\end{exam} 

\begin{theo} 
$V_{d,n} \subset \P^{n+1}$ is a small toric degeneration of a smooth 
Fano hypersurface $X_{d,n} \subset \P^n$ of degree $d$.
\label{sm-hyp}
\end{theo} 

\noindent
{\em Proof.} Let us first consider the case $n = 2d-2$. In this case 
the $2(d-1)$-dimensional fan $\Sigma_d$ defining the toric 
variety   $V_{d,2(d-1)}$ can be constructed as follows: 

Let $e_1, \ldots, e_{d-1}, f_1, \ldots, f_{d-1}$ be a $\Z$-basis of the 
lattice $\Z^{2(d-1)}$. We set  $e_{d} := -e_1 - \cdots - e_{d-1}$ 
and  $f_{d} := -e_1 - \cdots - f_{d-1}$. We denote by $h_{i,j}$ the 
sum $e_i + f_j$ ($i, j \in \{1, \ldots, n\})$. If $\Delta_d^*$ denotes 
the convex hull of $d^2$ points $h_{i,j}$, then the fan $\Sigma_d \subset 
N_{\R}$ consists of cones over faces of the reflexive polyhedron 
$\Delta_d^*$, where the integral lattice $N \subset \Z^{2(d-1)}$ 
is generated by all $d^2$ lattice vectors $h_{i,j}$ (the sublattice $N 
\subset \Z^{2(d-1)}$ coincides with $ \Z^{2(d-1)}$ unless $d =2$).

Using the combinatorial  characterisations of terminal toric 
singularities \cite{KMM}, one immediately obtains 
that all singularities of $V_{d,2(d-1)}$ are terminal, since 
the only $N$-lattice points on the faces of $\Delta_d^*$ are their vertices.  
If  $d \geq 3$, then the Picard group of 
$V_{d,2(d-1)}$ is generated by the class of the hyperplane section, i.e., 
$Pic(V_{d,2(d-1)})\cong  \Z$ and the anticanonical class of 
$V_{d,2(d-1)}$ is $d$-th multiple of the 
generator of $Pic(V_{d,2(d-1)})$. The latter can be show as follows: 

Consider a $(2d-3)$-dimensional face of $\Delta_d^*$ having vertices
\[ h_{i,j}, \;\; i \in \{1, \ldots, d-1\}, \; j \in \{1, \ldots, d\}. \] 
Then every $\Sigma_d$-piecewise linear function $\varphi\, : \, N_{\R} \to 
\R$, up to summing a linear function, can be normalized by 
the condition  
\[ \varphi(h_{i,j}) = 0, \;\; \forall i \in \{1, \ldots, d-1\}, \; 
\forall j \in \{1, \ldots, d\}. \] 
On the other hand, for any $j \neq j'$, $j, j' \in  \{1, \ldots, d\}$ 
four lattice points 
\[ h_{d,j}, h_{1,j}, h_{d,j'}, h_{1,j'} \]
generate a $3$-dimensional cone in $\Sigma_d$. Hence 
\[ \varphi(h_{d,j}) = \varphi(h_{d,j'}) \;\; \forall 
j, j' \in  \{1, \ldots, d\}. \]
This means that the space of all $\Sigma_d$-piecewise linear functions 
modulo linear functions is $1$-dimensional. The anticanonical 
class is represented by the $\Sigma_d$-piecewise linear function $\varphi_1$
taking values $1$ on each vector $h_{i,j}$ $i,j \in  \{1, \ldots, d\}$. 
Considering the difference 
\[ \varphi'_1 := \varphi_1 - \lambda, \]
where $\lambda$ is a linear function on $N_{\R}$ satisfying the conditions 
\[ \lambda(e_1) = \cdots =\lambda(e_{d-1}) =1, \; \lambda(e_d) = -(d-1), \; 
\lambda(f_1) = \cdots = \lambda(f_d) = 0, \]
we obtain a $\Sigma_d$-piecewise linear function having the properties
\[ \varphi_1'(h_{i,j}) = 0, \;\; \forall i \in \{1, \ldots, d-1\}, \; 
\forall j \in \{1, \ldots, d\} \]
and 
 \[ \varphi(h_{d,j}) = d\;\; \forall 
j \in  \{1, \ldots, d\}. \]      
So the class of $\varphi_1$ modulo linear functions is a $d$-th multiple
of a generator of $Pic(V_{d,2(d-1)})$.

The general case $n > 2(d-1)$ can be obtained by similar arguments using 
the fact that   $V_{d,n}$ is a projective cone over $V_{d,2(d-1)}$. 
In order to construct the required  
flat $1$-parameter family 
${\frak X}$ (cf.  \ref{def-small}), it suffices to consider a pencil of 
hypersurfaces of degree $d$ in  $\P^{n+1}$ joining  $X_{d,n}$ and 
 $V_{d,n}$. 
\hfill $\Box$ 

\begin{theo} 
Let $X_d \subset \P^{n+1}$ be a smooth Fano 
hypersurface of degree $d$. Then $X_d$ admits a small toric degeneration 
if and only if $n \geq 2d -2$. 
\label{hypersur}
\end{theo} 

\noindent
{\em Proof.}
By \ref{sm-hyp}, it suffices to show that $X_d$ does not 
admit a small toric degeneration   if $n< 2d-2$.  
Assume that $X_d$ admits a small toric degeneration $Y_d$. Then $Y_d$ 
is a toric hypersurface defined by a binomial equation $M_1 =M_2$ where 
$M_1$ and $M_2$ are monomials in $z_0, \ldots, z_{n+1}$ of degree $d$
$(z_0, \ldots, z_{n+1}$ are homogeneous coordinates on $\P^{n+1}$). 
If $n < 2d-2$, then at least one of the  monomials $M_1$ and  $M_2$ 
must be divisible by  $z_i^2$  for some $i \in \{0, \ldots, n+1\}$. 
We can assume that for instance $z_0^2$ divides $M_1$. If $z_k$ and 
$z_l$ are two variables appearing in $M_2$, then $n-2$-dimensional 
linear subspace 
\[ z_0 = z_k = z_l = 0 \]
is contained in $Sing(Y_d)$. This contradicts the fact that 
terminal singularities on $Y_d$ could appear only in codimension $\geq 3$ 
(see \ref{codim3}). 
\hfill $\Box$

Using \ref{codim3}, one  
immediately obtains:  

\begin{prop}
If $X$ is a smooth Del Pezzo surface, then $X$ admits a small toric
degeneration if and only if $X$ is itself a toric variety (i.e. 
$K_X^2 \geq 6$). 
\end{prop} 

As we have seen from \ref{sm-hyp}, a smooth quadric $3$-fold in 
$\P^4$ is an example of nontoric smooth Fano variety which 
admits a small toric degeneration. By \ref{hypersur}, cubic and 
quartic $3$-folds do not admit 
small toric degenerations. 
The compltete classification of smooth Fano $3$-folds has 
been obtained in \cite{C,I,MM1,MM2,MU}. It is natural to ask the following:

\begin{ques} 
Which $3$-dimensional nontoric smooth Fano varieties do admit small toric 
degenerations? 
\end{ques}

\section{The mirror construction} 

For our convenience, we assume $K= \C$. 

Let $X$ be a smooth Fano $n$-fold over $\C$ 
and $Y$ is its  small toric degeneration. 
The toric variety  $Y$ is defined by some  complete rational polyhedral fan 
 $\Sigma \subset N_{\R}$, where $N_{\R} = N \otimes \R$  is the 
real scalar extension of a $N \cong \Z^n$.
We denote by $Cl(Y)$ (resp. by $Pic(Y)$) the group of Weil (resp. 
Cartier) divisors on $Y$ modulo the rational equivalence. 
One has a canonical embedding 
$$\alpha\; : \; Pic(Y) \hookrightarrow  Cl(Y).$$
If   $\{ e_1, \ldots, e_k \} 
\subset N$ is the set of integral generators of $1$-dimensional cones 
in $\Sigma$, then $Cl(Y)$ is a finitely generated abelian 
group of rank $k-n$ and the convex hull of $
e_1, \ldots, e_k$ is a reflexive polyhedron $\Delta^*$ (for definition 
of reflexive polyhedra see \cite{BA}).
Assume that there exists a partition of the set 
$I = \{ e_1, \ldots, e_k \}$ into $r$ disjoint subsets 
$J_1, \ldots, J_r$ such that the union $D_i$ 
of toric strata in $Y$ corresponding 
to elements of $J_i$ is a semiample Cartier divisor on $Y$ for each 
$i \in \{1, \ldots, r\}$. Denote by $Z \subset Y$ a Calabi-Yau complete 
intersection of $r$ hypersurfaces $Z_i \subset Y$ defined by vanishing of 
generic global sections of ${\cal O}_Y(D_i)$.  
By \cite{BS} (see also \cite{Bo}), the mirrors $Z^*$ 
of Calabi-Yau complete intersections  $Z \subset Y$ are birationally 
isomorphic to affine complete intersections  in $(\C^*)^n = 
Spec\, \C[t_1^{\pm 1}, \ldots, 
t_n^{\pm 1}]$ defined by $r$  equations 
\[ 1 = \sum_{e_j \in J_i}^k  a_j{\bf t}^{e_j}, \;\; i \in \{1, \ldots, r\},\]
where $(a_1, \ldots, a_k) \in \C^k$ is  a general complex vector 
and ${\bf t}^{e_1}, \ldots, {\bf t}^{e_k}$ are Laurent monomials in 
variables $t_1, \ldots, t_n$ with the exponents 
$e_1, \ldots, e_k$.

\begin{dfn} 
{\rm  A complex vector  $(a_1, \ldots, a_k ) \in \C^k$ is called 
$\Sigma$-{\bf admissible},  if there exists a $\Sigma$-piecewise linear 
function 
\[ \varphi\; : \; N_{\R} \to \R, \]
(i.e., a continuous function such that $\varphi|_{\sigma}$ is linear 
for every $\sigma \in \Sigma$) having the property 
\[ \varphi(e_i) = \log|a_i|,\;\;\forall i \in \{1, \ldots, k\}. \] 
The set of all $\Sigma$-admissible vectors will be denoted 
by $A(\Sigma)$. 
} 
\end{dfn}

\begin{rem} 
{\rm It is easy to show that $A(\Sigma) \subset \C^k$ is 
an irreducible  closed subvariety which is isomorphic to an 
affine toric variety of dimension $rk\, Pic(Y) + n \leq k$.} 
\end{rem}

Now our generalization of the mirror construction 
from \cite{BS} to the case of Calabi-Yau complete intersections  
in a nontoric Fano variety $X$ can be formulated as follows: 
\medskip

\noindent
{\bf Generalized mirror construction:} 
{\em Mirrors $W^*$ of generic Calabi-Yau hypersurfaces $W \subset X$ 
are birationally isomorphic to the affine complete intersections
\[ 1 = \sum_{i=1}^k  a_i {\bf t}^{e_i}, \]
where ${\bf a}:= (a_1, \ldots, a_k)$ is a general 
point of $A(\Sigma)$.} 
\medskip  

\noindent
{\bf Monomial-divisor correspondence:} 
Let us explain the monomial-divisor mirror correspondence for this 
mirror construction (cf. \cite{AGM}). By \ref{def-small}(iii), the group 
$Pic(Y)$ can be canonically identified with $Pic(X)$. The image of the 
restriction homomorphism $Pic(X) \to Pic(W)$
defines a subgroup $G \subset Pic(W)$, whose elements correspond 
to monomial deformations of the complex structure on mirrors: 

{\em if $\psi$ is an integral $\Sigma$-piecewise linear function representing 
an element $\gamma \in G$, then 
the $1$-parameter family of hypersurfaces 
\[ 1 = \sum_{i=1}^k  t_0^{\varphi(e_i)} {\bf t}^{e_i},\;\; t_0 \in \C \]
defines the corresponding $1$-parameter deformation of the complex 
structure on $W^*$ via the deformation of the coefficients $a_i = 
 t_0^{\varphi(e_i)}$.} 
\medskip    

\noindent
{\bf The main period:} 
Let $R(\Sigma)$ the group of all vectors $(l_1, \ldots, l_k) \in 
{\Z}^k$ satisfying the condition 
$\sum_{i =1}^k l_i e_i = 0$
and 
$L(\Sigma) \subset R(\Sigma)$ be the semigroup consisiting of 
vectors $(l_1, \ldots, l_k) \in R(\Sigma)$ 
with nonnegative coordinates $l_i$ $(i =1, \ldots, k)$. 
There exists a canonical  
pairing 
$ \langle  *, * \rangle \; : \; R(\Sigma) \times Pic(Y) \to \Z $ 
which is the intersection pairing between $1$-dimensional cycles 
and Cartier divisors on $Y$. 
According to \cite{BS}, we can compute the main period in  the   
family of mirrors $W^*$ in our generalized mirror construction as 
follows
\[ \Phi_0({\bf a}) = \sum_{ {\bf l}= (l_1, \ldots, l_k) \in L(\Sigma)} 
\frac{ \langle  l, D_1 + \cdots + D_r \rangle!}{ \langle l, 
D_1\rangle! \cdots 
\langle l, D_r\rangle !} \prod_{i=1}^k a_i^{l_i}, \; \; {\bf a} \in 
A(\Sigma).  \] 
The condition ${\bf a} \in A(\Sigma)$ can be interpreted as 
a  specialization of  $GKZ$-hypergeometric series from \cite{BS}. 
\medskip

Some evidences in favor of our generalized  
mirror construction were presented 
in \cite{BCKS1,BCKS2}. For our  next  examples confirming  
the proposed generalized 
mirror construction we use the following simple combinatorial statement: 

\begin{prop} 
Let $S_d(m)$ be the set of all $d \times d$-matrices $K = (k_{ij})$ 
with nonnegative integral coefficients $k_{ij}$ satisfying 
the equations 
\[ \begin{pmatrix} 1 &  \cdots   & 1  \end{pmatrix} 
  \begin{pmatrix} k_{11} & \cdots & k_{1d} \\  
 \cdot & \cdots & \cdot  \\  \cdot & \cdots & \cdot  \\ 
\cdot & \cdots & \cdot  \\ k_{d1} & \cdots  & k_{dd} \end{pmatrix} =  
 \begin{pmatrix} m  & \cdots   & m  \end{pmatrix}  
\] 
and 
\[  
  \begin{pmatrix} k_{11} & \cdots & k_{1d} \\  
 \cdot & \cdots & \cdot  \\  \cdot & \cdots & \cdot  \\ 
\cdot & \cdots & \cdot  \\ k_{d1} & \cdots  & k_{dd} \end{pmatrix} 
\begin{pmatrix} 1 \\   \cdot \\ \cdot \\ \cdot  \\  1  \end{pmatrix}=  
 \begin{pmatrix} m  \\ \cdot \\   \cdot 
\\ \cdot  \\  m  \end{pmatrix}.  
\] 
Then 
\[ \sum_{K \in S_d(m)} \frac{(m!)^d}{\prod_{i,j =1}^{d} (k_{ij})!} = 
\frac{(dm)!}{(m!)^d}. \] 
\label{comb-f}
\end{prop}

\noindent
{\em Proof.} Let $A$ be the set $\{1, 2, \ldots, dm \}$ of first $dm$ 
natural numbers. We fix a 
splitting  $A$ into the  disjoint union of $d$ subsets
\[ A_i := \{ (i-1)m +1, (i-1)m +2, \ldots, im \}, \;\; i =1, \ldots, d \]
consising of $m$ elements. 
Let $\beta : 
A = B_1 \cup \cdots \cup B_d$ be an  arbitrary representation of $A$ 
as a disjoint union of the subsets $B_1, \ldots, B_d$ with the property 
$|B_1| = \cdots = |B_d| =m$. Then every  such a representation  defines a 
matrix $K(\beta) =(k_{ij}(\beta)) \in S_d(m)$ as 
follows: 
\[ k_{ij}(\beta) := |A_i \cap B_j|, \;\; i, j \in \{1, \ldots, d\}. \]
For a fixed matrix $K \in S_d(m)$ there exist exactly 
\[  \prod_{j=1}^d \frac{(m!)}{\prod_{i =1}^{d} (k_{ij})!} \]
ways to construct a representation $\beta$ of $A$ as a dusjoint union of 
$m$-element subsets $B_1, \ldots, B_d$ such that $K = K(\beta)$. 
Therefore,
\[ \sum_{K \in S_d(m)} \frac{(m!)^d}{\prod_{i,j =1}^{d} (k_{ij})!} \]
is the total number of ways to split $A$ into a disjoint union of 
$m$-element subsets $B_1, \ldots, B_d$. On the other hand, this number
is equal to the multinomial
\[ \frac{(dm)!}{(m!)^d} .\]
\hfill $\Box$ 

\begin{exam} 
{\em Let $W$ be a generic Calabi-Yau complete intersection of 
two hypersurfaces $V_d, V_d'$ in $\P^{2d-1}$. 
By \ref{sm-hyp}, we can construct a small toric degeneration 
of one smooth hypersurface  $V_d'$ to the 
$2(d-1)$-dimensional toric variety  $Y_d \subset \P^{2d-1}$ 
\[ z_0z_1 \cdots z_{d-1} = z_d z_{d+1} \cdots z_{2d-1}.  \]
Using an explicit description of the Picard group $Pic(Y_d)$ from 
the proof of \ref{sm-hyp},  
our generalized mirror construction suggests that mirrors $W^*$ 
for $W$ are birationally isomorphic to the  
affine hypersurfaces $Z_F$ in the algebraic 
torus 
$$Spec\, \C[ t_1^{\pm1},\ldots, t_{d-1}^{\pm1}, 
u_1^{\pm1},\ldots, u_{d-1}^{\pm1}]$$
defined by the $1$-parameter family of the equations
\[ 1 = F(t_1,\ldots, t_{d-1},u_1, \ldots, u_{d-1},z) =  
\sum_{i=1}^{d-1} \sum_{j=1}^{d-1}t_iu_j +  (u_1 \cdots u_{d-1})^{-1} 
\left( \sum_{i=1}^{d-1} t_i \right)  + \]
\[ + z(t_1 \cdots t_{d-1})^{-1} 
\left(  (u_1 \cdots u_{d-1})^{-1} + 
\sum_{i=1}^{d-1} u_j \right) ,\; \;\;\;  z \in \C \] 

On the other hand, it is known via a toric  mirror construction 
for Calabi-Yau complete intersection $W = V_d \cap V_d'$ 
(see \cite{BS}) that the power series 
\[ \Phi_0(z) = \sum_{m \geq 0} \frac{(dm!)^2}{(m!)^{2d}} z^m \] 
generates the Picard-Fuchs $D$-module discribing the  
quantum differential system. 
Now we compare our generalized mirror construction with the 
known one from \cite{BS} computing the main period of the family $Z_F$ 
by the Cauchy residue formula:   
\[ \Psi_F(z) := \frac{1}{(2\pi\sqrt{-1})^{2(d-1)}} \int_{\Gamma} \frac{1}{1 -  
F({\bf t}, {\bf u}, z)} \frac{{\bf dt}}{{\bf t}} 
\wedge \frac{{\bf du}}{\bf u} = 
1 + a_1z + a_2z^2 + \cdots, \]
\[  \frac{{\bf dt}}{{\bf t}}:=   \frac{dt_1}{t_1} \wedge \cdots   
\wedge \frac{dt_{d-1}}{t_{d-1}}, \; \; \; 
\frac{{\bf du}}{{\bf u}}:=   \frac{du_1}{u_1} \wedge \cdots   
\wedge \frac{du_{d-1}}{u_{d-1}}, \]
where the coefficients $a_m$ of the power series $\Psi_F(z)$ 
can be computed by the formula 
\[ a_m = \sum_{K \in S_d(m)} \frac{(dm)!}{\prod_{i,j =1}^{d} (k_{ij})!}. \]
Using \ref{comb-f}, we obtain 
that 
\[ a_m =  \frac{(dm!)^2}{(m!)^{2d}}, \]
i.e., the power series $\Psi_F(z)$ coincides with $\Phi_0(z)$ and therefore
our generalized mirror construction agrees with the already known one 
from \cite{BS}. 
 
For the special case $d=3$, we obtain a description 
for  mirrors $W^*$ of  complete intersections $W$ of two cubics in $\P^5$  
as   smooth compactifications 
of  hypersurfaces in the $4$-dimensional algebraic 
torus $$Spec\, \C[ t_1^{\pm1},t_2^{\pm1}, u_1^{\pm1},u_2^{\pm1}]$$
defined by the $1$-parameter family of the equations
\[ 1 = F(t_1,t_2,u_1,u_2, \lambda) =  t_1u_1 + t_1u_2 + t_1(u_1u_2)^{-1} + 
t_2u_1 + t_2u_2 +  t_2(u_1u_2)^{-1} + \]
\[ + z(t_1t_2)^{-1}(u_1 + u_2 + 
(u_1u_2)^{-1}),\; \;\;\;  z \in \C.  \]
This discription of mirrors is different from the one proposed by Libgober and 
Teilelbaum in \cite{LT}, but it seems that both constructions are equivalent 
to each other.}
\end{exam} 

Now we want to suggest some problem which naturally arise from 
the  proposed generalized mirror construction.

\begin{prob} 
Check the topological mirror duality test
\[ E_{\rm st}(W^*; u,v) = 
(-u)^n  E_{\rm st}(W; u^{-1},v) \]
for the above generalized mirror construction. Here $E_{\rm st}$ is 
the stringy $E$-function introduced in \cite{B1}.  
\end{prob}

\begin{rem} 
{\rm The main difficulty of this checking arises from the fact that 
the affine complete intersections in the above  mirror construction 
are not {\em generic}. For $\Delta^*$-regular affine hypersurfaces 
there exists explicit combinatorial formula for their  $E$-polynomials 
(see \cite{BB1}). However, the affine hypersurfaces in our mirror 
construction are not $\Delta^*$-regular and no explicit formula for their  
$E$-polynomials (or Hodge-Deligne numbers) is known so far. 
} 
\end{rem}

\begin{prob} 
Generalize the method of Givental \cite{G1,G2,G3} for computing  
Gromov-Witten invariants of complete intersections in 
smooth Fano varieties $X$ admitting small toric degenerations.
\end{prob} 

\begin{rem} 
{\rm If $X$ is a smooth Fano $n$-fold admitting a small toric degeneration 
$Y$, then one can not expect that there exists a $\C^*$-action on $X$. So 
the equivariant arguments from  \cite{G1} can not be applied directly 
to $X$. However, 
one could try to use equivariant Gromov-Witten theory for the 
ambient projective space $\P^m$ containing both $X$ and $Y$  and to 
show that the virtual fundamental classes corresponding to 
$Y$ and $X$ are the same. It seems that small quantum cohomology of $Y$ 
carry complete information about the subring in the small quantum cohomology 
ring $QH^*(X)$ generated by the classes of divisors. This would give 
an explicit description of such a subring (see \cite{ST}) as well as 
of its  gravitational version via Lax operators (see \cite{EHX}).  
}  
\end{rem}

\end{document}